\documentclass{amsart}

\usepackage{amsmath,amssymb}
\usepackage[all]{xy}
\usepackage{hyperref}

\newtheorem{theorem}{Theorem}[section]
\newtheorem{lemma}[theorem]{Lemma}
\newtheorem{proposition}[theorem]{Proposition}

\theoremstyle{definition}
\newtheorem{definition}[theorem]{Definition}
\newtheorem{example}[theorem]{Example}

\theoremstyle{remark}
\newtheorem{remark}[theorem]{Remark}
\newtheorem{notation}[theorem]{Notation}

\newcommand {\R} {\ensuremath {\mathbb{R}} }
\newcommand {\isom} {\ensuremath {\cong} }
\newcommand {\tensor} {\ensuremath {\otimes} }
\newcommand {\incl} {\ensuremath {\hookrightarrow} }
\newcommand {\isomto} {\ensuremath {\xrightarrow{\isom}} }
\newcommand {\xto}[1] {\ensuremath {\xrightarrow{#1}} }
\newcommand {\Top} {\ensuremath {\mathbf{Spaces}} }
\newcommand {\Set} {\ensuremath {\mathbf{Set}} }
\newcommand {\sSet} {\ensuremath {\mathbf{sSet}} }
\newcommand {\cSet} {\ensuremath {\mathbf{cSet}} }
\newcommand {\LPS} {\ensuremath {\mathbf{LPS}} }
\newcommand {\ALPS} {\ensuremath {\mathbf{A \downarrow \LPS}} }
\newcommand {\yAsPreLPS} {\ensuremath {\mathbf{\bar{y}(A)\downarrow \sPre(\LPS)}} }
\newcommand {\Pospc} {\ensuremath {\mathbf{PoSpaces}} }
\newcommand {\APospc} {\ensuremath {\mathbf{A \downarrow PoSpaces}} }
\newcommand {\Sh}  {\ensuremath {\mathbf{Shv} } }
\newcommand {\sPre}  {\ensuremath {\mathbf{sPre} } }
\newcommand {\sPreLPS} {\ensuremath {\mathbf{\sPre(\LPS)}} }
\newcommand {\ShLPS} {\ensuremath {\mathbf{\Sh(\LPS)}} }
\newcommand {\C} {\ensuremath {\mathbb{C}} }
\newcommand {\ShC} {\ensuremath {\Sh(\C)} }
\newcommand {\PreC} {\ensuremath {\Set^{\C^{op}}} }
\newcommand {\sPreC} {\ensuremath {\sSetCop}}
\newcommand {\OZ} {\ensuremath {\mathcal{O}(Z)}}
\newcommand {\EtaleZ} {\ensuremath {\mathbf{Etale}(Z)}}
\newcommand {\ShOZ} {\ensuremath{\Sh(\OZ)}}
\newcommand {\mmBook} {maclaneMoerdijkBook1992}
\newcommand {\SetCop} {\ensuremath {\Set^{\C^{\op}}} }
\newcommand {\SetLPSop} {\ensuremath {\Set^{\LPS^{\op}}} }
\newcommand {\sSetCop} {\ensuremath {\sSet^{\C^{\op}}} }
\newcommand {\E} {\ensuremath {\mathcal{E}} }
\newcommand {\Iff} { if and only if }
\newcommand {\Rmphi} {\bar{M}_{\phi}}
\newcommand {\Rnphi} {\bar{N}_{\phi}}
\newcommand {\adjn} {\leftrightarrows}
\newcommand {\opensubset} {\stackrel{\subset}{{\scriptscriptstyle \open}}}
\newcommand {\dI} {\vec{I}}
\newcommand {\from} {\leftarrow}
\newcommand{\coatl}[1]{\mathcal{#1}}

\newcommand {\M} {\ensuremath {\mathbf{\mathcal{M}}}}
\newcommand {\N} {\mathbf{\mathcal{N}}}
\newcommand {\AM} {\ensuremath {\mathbf{A \downarrow \mathcal{M}}}}

\DeclareMathOperator{\Mor}{Mor}
\DeclareMathOperator{\op}{op}
\DeclareMathOperator*{\colim}{colim}
\DeclareMathOperator{\cod}{cod}
\DeclareMathOperator{\dom}{dom}
\DeclareMathOperator{\Id}{Id}
\DeclareMathOperator{\im}{im}
\DeclareMathOperator{\germ}{germ}
\DeclareMathOperator{\open}{open}

\newdir{ >}{{}*!/-10pt/\dir{>}}

\newcommand{\es}{\mathcal{E}}
\newcommand{\fs}{\mathcal{F}}
\newcommand{\os}{\mathcal{O}}
\newcommand{\deq}{\stackrel{def}{=}}

\newcommand{\pre}[1]{\set^{#1^{op}}}
\newcommand{\sset}{\sSet}
\newcommand{\spre}[1]{\sset^{#1^{op}}}
\newcommand{\stal}[1]{stalk_{#1}}
\newcommand{\set}{\Set}

\newcommand{\leirom}{\renewcommand{\labelenumi}{\textit{(\roman{enumi})}}}
\newcommand{\leiarab}{\renewcommand{\labelenumi}{(\arabic{enumi})}}

\newenvironment{enumerateroman}
{\leirom \begin{enumerate}}
{\end{enumerate} \leiarab}

\begin{document}

\title{A model category for local po-spaces}

\author{Peter Bubenik}
\email{p.bubenik@csuohio.edu}
\address{  Department of Mathematics\\
  Cleveland State University\\
  2121 Euclid Ave. RT 1515\\
  Cleveland OH, 44115-221\\
  USA}
\thanks{This research was partially funded by the Swiss National Science
  Foundation grant 200020-105383.}

\author{Krzysztof Worytkiewicz}
\email{kworytki@uwo.ca}
\address{Department of Mathematics\\
  University of Western Ontario\\
  Middlesex College\\
  London, Ontario N6A 5B7\\
  Canada}



\begin{abstract}
Locally partial-ordered spaces (local po-spaces) have been used to
model concurrent systems. 
We provide equivalences for these spaces by constructing a model
category containing the category of local po-spaces.
We show the category of simplicial presheaves on local po-spaces can
be given Jardine's model structure, in which we identify the weak
equivalences between local po-spaces.
In the process we give an equivalence between the category of sheaves
on a local po-space and the category of {\'e}tale bundles over a local
po-space. 
Finally we describe a localization that should
provide a good framework for studying concurrent systems.
\end{abstract}

\maketitle

\section{Introduction} \label{sectionIntro}

The motivation for this paper stems from the study of concurrent 
processes accessing
shared resources. 
Such systems were originally described by discrete 
models based on graphs, possibly equipped with some additional 
information \cite{calcs}. The precision of these models 
suffers, however, from an inaccuracy in distinguishing between concurrent and 
non-deterministic executions. It turned out that a satisfactory way to organize
this information can be based on cubical sets, giving rise to the notion of 
\emph{Higher-Dimensional Automata} or HDA's 
\cite{Goubault:1996:DTT,eric:cmcim}. HDA's live in 
slice categories of $\cSet$, the category of cubical sets and their morphisms.

A different view, which has its origins in Dijkstra's notion of 
\emph{progress graphs} \cite{dijk68}, takes the flow of time into account. 
The difficulty
here is to adequately model the fact that time is irreversible as far as
computation is concerned. On the other hand, one would like to identify
execution paths corresponding to (at least) the same sequence of 
acquisitions of shared 
resources. However, in order not to lose precision, this notion of homotopy 
is also subject to the constraint above of the irreversibility of time.
There are two distinct approaches, both based on topological spaces. 

One approach, advocated by P. Gaucher, 
is to topologize the sets of paths between the states of an automaton, which
technically amounts to an enrichment with no units \cite{gau1}. 
The intuition behind the
setup is to distinguish between \emph{spatial} and 
\emph{temporal} deformations of computational 
paths. The related framework of \emph{Flows} has clear technical advantages 
from a (model-)categorical point of view.

The other approach, advocated by Fajstrup, Goubault, Raussen and 
others, is to
topologize partially ordered states of automata. Such objects
are called partially-ordered spaces or \emph{po-spaces} (also
\emph{pospaces})\footnote{M. Grandis uses a related
  approach~\cite{grandisDHTi} in which the underlying topological
  space comes with a class of directed paths. However these spaces are
  not partially-ordered, even locally.}. 
The advantage of using po-spaces is that there 
is a very simple
and intuitive way to express directed homotopy or \emph{dihomotopy}
\cite{goubaultSomeGPiCT,fgrAlgebraicTaCpreprint}.

However, the
price paid is that po-spaces cannot model executions of (concurrent) programs
with loops. The solution is to order the underlying topological space
only \emph{locally}. Such objects are called
\emph{local po-spaces} and the notion of dihomotopy becomes more 
intricate in this context. Nevertheless, practical reasons like 
tractability
call for a good notion 
of equivalence in the category of
local po-spaces. Put differently, it would be useful to be able to
replace a given
local po-space model with a simpler local po-space which nevertheless
preserves the relevant computer-scientific properties. 

In this paper, we study these questions in the  framework of Quillen's
(closed) model categories
\cite{quillenHomotopicalAlgebra,hoveyBook,hirschhornBook}.  
Briefly, a model category is a category with all small limits and
colimits and three distinguished classes of morphisms called \emph{weak
  equivalences}, \emph{cofibrations} and \emph{fibrations}.
Weak equivalences that are also cofibrations or fibrations are called
\emph{trivial cofibrations} and \emph{trivial fibrations}, respectively.
These morphisms satisfy four axioms that allow one to apply the
machinery of homotopy theory to the category. This machinery allows a
rigorous study of equivalences.
We remark that there are other frameworks for studying equivalence.
However model categories have the most developed theory, and have
succeeded in illuminating many diverse subjects.

Our aim is to construct a model category of locally partial-ordered
spaces as a foundation for the study of concurrent systems.
This is technically difficult because locally partial-ordered spaces
are not closed under taking colimits.
We will define a category \LPS of local po-spaces, which embeds into
the category \sPreLPS of simplicial presheaves on local po-spaces.
The objects of \sPreLPS are contravariant functors from \LPS to the
category of simplicial sets and the morphisms are the natural
transformations. 
This embedding is given by a Yoneda embedding (see
Definition~\ref{defnSimplicialYoneda}),
\[
\bar{y}: \LPS \to \sPreLPS.
\]

We now briefly describe some technical conditions on model categories
which strengthen our theorems. For more details
see Definitions \ref{defnCofGen} and \ref{defnCellular} and
\cite{hoveyBook,hirschhornBook}. 
A model category is \emph{proper} if the weak equivalences are closed under
both pushouts with cofibrations and pullbacks with fibrations.
It is \emph{left proper} if the first condition holds.
A model category is \emph{cofibrantly generated} if the model category
structure is induced by a set of generating cofibrations and a set of
generating trivial cofibrations, both of which permit the small object
argument.
A \emph{cellular} model category is a cofibrantly generated model
category in which the cell complexes are well behaved. 
A \emph{simplicial} model category $\mathcal{M}$ is a model category
enriched over simplicial sets, which for any $X \in \mathcal{M}$ and any
simplicial set $K$ has objects $X \tensor K$ and $X^K$ which satisfy
various compatibility conditions.

\begin{theorem} \label{thmMain}
The category \sPreLPS has a proper, cellular, simplicial
model structure in which  
\begin{itemize}
\item the cofibrations are the monomorphisms, 
\item the weak equivalences are the \emph{stalkwise equivalences}, and 
\item the fibrations
are the morphisms which have the right lifting property with respect
to all trivial cofibrations.
\end{itemize}
Furthermore among morphisms coming from \LPS (using the Yoneda
embedding $\LPS \incl \sPreLPS$), the weak equivalences are precisely the
isomorphisms.
\end{theorem}

The model structure on \sPreLPS is Jardine's model
structure~\cite{jardineSimplicialPresheaves,jardineBooleanLocalization}
on the category of simplicial presheaves on a small \emph{Grothendieck
  site}.  
We show that $\ShLPS$ is a \emph{Grothendieck topos} which has
\emph{enough points}.  
Under this condition, Jardine showed that the weak equivalences are
the \emph{stalkwise equivalences}.

This model category can be thought of as a localization of the
universal injective model category of local
po-spaces~\cite{joyal:letter,duggerUHT,duggerHollanderIsaksenHypercoversaSP}.
While in general the weak equivalences are interesting and
nontrivial~\cite{jardineSimplicialPresheaves}, this is not true for
those coming from $\LPS$.
To obtain a more interesting category from the point of view of
concurrency we would like to localize with respect to directed
homotopy equivalences.
In~\cite{bubenikContextBRICS} it is argued that the relevant equivalences
are the directed homotopy equivalences relative to some
\emph{context}. 
The context is a local po-space $A$ and the directed homotopy
equivalences rel $A$ are a set of morphisms in $\ALPS$.

We combine this approach with Theorem~\ref{thmMain} as follows.
First we remark that $A$ embeds in $\sPreLPS$ as $\bar{y}(A)$.
Next the model structure on $\sPreLPS$ induces a model structure on
the coslice category $\yAsPreLPS$.
Finally one can take the left Bousfield localization of this model
category with respect to the directed homotopy equivalences rel $A$.

\begin{theorem} \label{thmLocalized}
Let $\mathcal{I} = \{ \bar{y}(f) \ | \ f \text{ is a  directed
    homotopy equivalence rel }A \}$.
Then the category \yAsPreLPS has a left proper, cellular
model structure in which  
\begin{itemize}
\item the cofibrations are the monomorphisms, 
\item the weak equivalences are the $\mathcal{I}$-local equivalences, and
\item the fibrations are those morphisms which have the right lifting
  property with respect to monomorphisms which are $\mathcal{I}$-local
  equivalences. 
\end{itemize}
\end{theorem}

Recall that, given a topological space $Z$, \emph{{\'e}tale bundles over $Z$}
are maps $W\to Z$ which are local
homeomorphisms. Let $\OZ$ be $Z$'s locale of open subsets and recall
that sheaves over $Z$ are functors $\OZ^{op} \rightarrow \Set$ that
enjoy a good gluing property. There is a well-known  
correspondence between \'etale bundles and sheaves.  
We establish a
directed version of this correspondence, which may be of independent interest.

\begin{theorem} \label{thmEtaleSheafCorrespondence}
Let $Z \in \LPS$.
Let $\EtaleZ$ be the category of di-\'etale bundles over $Z$, i.e. 
the category of bundles which are local dihomeomorphisms. 
Let $\OZ$ be the category of open subobjects of $Z$. 
There is an equivalence of categories:
\[
\Gamma: \EtaleZ \adjn \ShOZ: \Lambda.
\]
\end{theorem}

\noindent
\textbf{Acknowledgments.} The authors would like to thank Eric Goubault,
Emmanuel Haucourt, Kathryn Hess and Phil Hirschhorn for helpful
discussions and suggestions.

\tableofcontents

\section{Background}

This section contains some known definitions and facts we build on. We start
by
stating the definition of a model category in
 subsection~\ref{subsec-modelCat}.
Next we review the basics on presheaves in subsection
\ref{subsec-presheaves} and on sheaves in subsection \ref{subsec-sheaves}.
We then recall the notions of topoi and geometric morphisms in subsection
\ref{subsec-topoi} and of stalks in subsection \ref{subsec-stalks}. Our main 
reference for this material is \cite{\mmBook}. 
Subsection \ref{subsec-spresheaves} is devoted to some important model
structures on $\spre{\C}$, the category of simplicial presheaves over a 
category $\C$. The material is drawn from \cite{jardineSimplicialPresheaves,jardineBooleanLocalization,duggerHollanderIsaksenHypercoversaSP}.

\subsection{Model categories} \label{subsec-modelCat}

Recall that a morphism   $i:A \to B$ has the left lifting property
with respect to a morphism $p:X \to Y$ if in every commutative diagram
\begin{center}$
\xymatrix{
A \ar[r] \ar[d]_i & X \ar[d]^p\\
B \ar[r] & Y}
$\end{center}
there is a morphism $h:B \to X$ making the diagram commute. 
Also $f$ is a retract of $g$ if there is a commutative diagram:
\[
\xymatrix@=0.3cm{
A \ar[dd]_f \ar@{=}[rr] \ar[dr] & & A \ar[dd]^f \\
& X \ar[dd]_(0.3){g} \ar[ur] \\
B \ar@{=}'[r][rr] \ar[dr] & & B \\
& Y \ar[ur]
}
\]
\begin{definition}
A \emph{model category} is a category with all small
limits and colimits that has three distinguished classes of
morphisms:  $\mathcal{W}$, called the \emph{weak equivalences};
$\mathcal{C}$, called the \emph{cofibrations}; and $\mathcal{F}$, called
the \emph{fibrations}, which together
satisfy the axioms below. We remark that morphisms in $\mathcal{W}$ $\cap$
$\mathcal{C}$, and $\mathcal{W}$
$\cap$ $\mathcal{F}$, are called trivial cofibrations and trivial fibrations,
respectively. 
\begin{enumerate}
\item Given composable morphisms $f$ and $g$ if any of the two
  morphisms $f$, $g$, and $g \circ f$ are in $\mathcal{W}$, then so is the third.
\item If $f$ is a retract of $g$ and $g$ is in $\mathcal{W}$,
  $\mathcal{C}$ or $\mathcal{F}$, then so is $f$.
\item Cofibrations have the left-lifting property with respect to trivial fibrations, and trivial cofibrations have the
 left-lifting property with respect to fibrations.
\item Every morphism can be factored as a cofibration followed by a
  trivial fibration, and as a trivial cofibration followed by a fibration. These factorizations are functorial.
\end{enumerate}
\end{definition}

\subsection{Presheaves}\label{subsec-presheaves}

Recall that a presheaf $P$ on $\C$  is just 
a functor 
$P \in \PreC$. In particular, ``hom-ing''
\[
\begin{array}{llcl}
\C(\_,C): & \C^{op} & \rightarrow & \Set
\\
& X & \mapsto & \C(X,C)
\end{array}
\]
gives rise to a presheaf and further to the Yoneda embedding
\[
\begin{array}{llcl}
y: & \C & \rightarrowtail & \PreC
\\
& C & \mapsto & \C(\_,C).
\end{array}
\]
This embedding is \emph{dense}, i.e. 
\[
P \cong \colim (y \circ \pi)
\]
canonically for any 
presheaf $P$, where $\pi: (y \downarrow P) \rightarrow \C $ is the 
projection from the comma-category $y \downarrow P$. Recall that a presheaf 
in the image of the 
Yoneda-embedding (up to equivalence) is called \emph{representable}.

\subsection{Sheaves}\label{subsec-sheaves}

\begin{definition}\label{def-groth}
A \emph{sieve} on $M \in \C$ is a subfunctor $S \subseteq \C(\_,M)$.
A \emph{Grothendieck topology} $J$ on $\C$ assigns to each $M \in \C$ a 
collection $J(M)$ of sieves on $M$
such that
\begin{enumerateroman}
\item (maximal sieve) $\C(\_,M) \in J(M)$ for all $M \in \C$;
\item (stability under pullback) if $g: M \to N$ and $S \in J(N)$, then 
$(g\circ \_)^*(S) \in J(M)$ as given by
\begin{center}
$
\xy
\xymatrix{
(g \circ \_)^* (S) \ar[r] \ar@{ >->}[d] & S \ar@{ >->}[d] 
\\
\C(\_,M) \ar[r]_{(g \circ \_)} & \C(\_,N)
}
\POS( 5.5 , -6);\POS( 5.5,-3.6 ) \connect@{-}
\POS( 5.5 , -6);\POS( 3.1,-6 ) \connect@{-}
\endxy
$
\end{center}
\item (transitivity) if $S \in J(M)$ and $R$ is a sieve on $M$ such that 
$(f\circ \_)^*(R) \in J(U)$ for all $f:U \to M$ in the image of $S$, then $R \in J(M)$;
\end{enumerateroman}
\end{definition}

We say that a sieve $S$ on $M$ is a \emph{covering sieve} or a \emph{cover of}
$M$ whenever $S \in J(M)$.

\begin{remark}
Unwinding definition \ref{def-groth} pinpoints a sieve as a 
right ideal, i.e. a 
set of arrows $S$ with 
codomain $M$ such that $f \in S \implies
f \circ h \in S$ whenever the codomain of $h$, $\cod(h) = \dom(f)$,
the domain of $f$. From this
point of view, pulling back a sieve $S$ on $M$ by an
arrow $N\xto{f} M$
amounts to building the set
\[
f^*(S) \deq \{h | \ \cod(h) = N, \ f \circ h \in S\}.
\]
It is then immediate how to rephrase a Grothendieck topology in terms of right
ideals. 
\end{remark}

\begin{definition}
Let $J$ be a Grothendieck topology on $\C$. A presheaf $P \in \PreC$ is a \emph{sheaf} with respect to $J$ provided any natural transformation
$\theta: S \Rightarrow P$ uniquely extends through $y(M)$ as in
\begin{center}
$\xymatrix{
S \ar[r]^\theta \ar@{ >->}[d] & P
\\
y(M) \ar@{.>}[ur]
}$
\end{center}
for all $S \in J(M)$ and all $M \in \C$. $J$ is \emph{subcanonical} if
the representable presheaves are sheaves.
\end{definition}

\begin{remark} \label{defnMatchingFamily}
Let $\theta: S \rightarrow P$ be a natural transformation from a sieve $S$ 
to a presheaf $P$. If one sees $S$ as a right ideal
$S = \{u_j: M_j \to M\}$, then $\theta$ amounts 
to a
function that assigns to every $u_j:M_j \to M \in S$ an element
$a_j \in P(M_j)$ such that
\[
P(v)(a_j) = a_k
\]
for all $v: M_k \rightarrow M_j$ and for all $u_k = u_j \circ v \in S$. 
Such a function is called 
a \emph{matching family} for $S$ of elements of $P$. A matching family 
$a_j \in P(M_j)$ admits an \emph{amalgamation}
$a \in P(M)$ 
if 
\[
P(u_j)(a) = a_j
\]
for all $u_j: M_j \to M \in S$. 
From this point of view, the Yoneda lemma characterizes a \emph{sheaf} as
a presheaf such 
that every matching family has
a unique amalgamation for all $S \in J(M)$ and all $M\in \C$.
\end{remark}

A Grothendieck topology is a huge object. In practice, a generating
device is used.

\begin{definition} \label{defnBasisGrTop}
A \emph{basis} $K$ for a Grothendieck topology assigns to each 
object $M$ a collection $K(M)$
of families of morphisms with codomain $M$ such that
\begin{enumerateroman}
\item all isomorphisms $f: U \to M$ are contained in $K(M)$,
\item \label{itemStabilityAxiom}
given a morphism $g: N \to M \in \C$ and $\{f_i: U_i \to M\} \in K(M)$,
then the family of pullbacks $\{\pi_2: U_i \times_M N \to N\} \in
K(N)$, and
\item
given $\{f_i: U_i \to M\} \in K(M)$ and for each $i$, $\{h_{ij}:
A_{ij} \to U_i\} \in K(U_j)$, then the family of composites $\{f_i \circ
h_{ij}: A_{ij} \to M\} \in K(M)$.
\end{enumerateroman}
\end{definition}

\begin{remark}
Given a basis $K$ for a Grothendieck topology one generates the
corresponding Grothendieck topology $J$ by defining
\[ V \in J(M) \iff \text{there is } U \in K(M) \text{ such that } U
\subset V.
\]
As expected, the sheaf condition can be rephrased in terms of a basis.
\end{remark}

As an example, consider  the case $\C = \os(X)$ with
$X$ a topological space and $\os(X)$ its locale of opens. The basis of the
\emph{open-cover} (Grothendieck) topology is, as expected, given by
open coverings of the opens.

\begin{theorem}\label{theo-sheafi}
Let $\Sh(\C,J)$ be the full subcategory of
$\PreC$ whose objects are sheaves for $J$.
The inclusion functor $i: \Sh(\C,J) \to \SetCop$ has a left adjoint $a$
called the \emph{associated sheaf functor} or \emph{sheafification}. This
left adjoint preserves finite limits.
\end{theorem}

Theorem \ref{theo-sheafi} is listed as Theorem III.5.1 in \cite{\mmBook}. There
are several equivalent ways to construct the associated sheaf functor, the
most classical one being the ``plus-construction'' applied twice. 

\begin{remark}
A cover on $M$ amounts to a cocone in $\C$ with vertex $M$. The associated
sheaf functor maps these cocones onto colimiting ones. Moreover, it is 
universal
with respect to this property. 
\end{remark}

\subsection{Topoi}\label{subsec-topoi}

\begin{definition}
A category $\E$ \emph{has exponentials} provided that for all $X \in
\E$, the functor $\_ \times X: \E \to \E$ has a right adjoint denoted
$(\_)^X$, so that 
\[
\E(Y \times X, Z) \isom \E(Y, Z^X).
\]
Suppose now $\E$ has a terminal object $1$, and has finite limits. 
A \emph{subobject classifier} is a monomorphism
$\text{true}: 1 \rightarrowtail \Omega$ such that for every 
monomorphism $s: S \rightarrowtail X$, there 
is a unique
morphism $\phi_S$ such that pullback of $\text{true}$ along $\phi_S$ yields
$s$:
\begin{center}
$
\xy
\xymatrix{
S \ar[r] \ar@{ >->}[d]_s & 1 \ar@{ >->}[d]^{\text{true}} 
\\
X \ar[r]_{\phi_S } & \Omega
}
\POS( 4 , -4);\POS( 4,-1.6 ) \connect@{-}
\POS( 4 , -4);\POS( 1.6,-4 ) \connect@{-}
\endxy
$
\end{center}
The category $\E$ is a \emph{topos} if it has exponentials and a 
subobject classifier.
\end{definition}

A subobject classifier is obviously unique
(up to isomorphism).
Furthermore, a topos has all finite colimits, though this is not easy
to prove.
It would take pages to enumerate all the remarkable
features of a topos, see \cite{JohnstonePT:topt} for an introduction to 
the lore of the material. Let us just say that topoi as introduced by 
Grothendieck and his collaborators had a very strong
algebro-geometrical flavor \cite{sga4}, yet the rich structure is
relevant not only for for algebraic geometers but for logicians as well
\cite{LawvereFW:diss,LawvereFW:eletcs,Law73}.

\begin{definition}
A site $(\C,J)$ is a small category $\C$ equipped with a Grothendieck 
topology $J$. A \emph{Grothendieck topos} is a category equivalent 
to the
category $\Sh(\C,J)$ of sheaves on $(\C,J)$.
\end{definition}

The following are well known.

\begin{proposition}
\begin{enumerate}
\item
A Grothendieck topos is a topos;
\item
$\set$ is a topos;
\item
$\SetCop$ is a topos for any $\C$.
\end{enumerate}
\end{proposition}

\begin{definition}
Let $\es$ and $\fs$ be topoi. A \emph{geometric morphism} 
$g: \fs \rightarrow \es$ is a pair of adjoint functors
\[
\xymatrix{
\es \ar@<1ex>[r]^{g^*}
     & \fs \ar@<1ex>[l]^{g_*}_{\perp} }
\]
such
that the left adjoint $g^*$ is left-exact (that is, it preserves finite
limits). The right adjoint is called \emph{direct image} and the left
one \emph{inverse image}.  
\end{definition}

As an example, $i: \Sh(\C,J) \hookrightarrow \PreC$ is the direct image part
of a geometric morphism. Notice that the convention for a geometric 
morphism is to have the direction
of its direct image part. 

\begin{definition}
A (geometric) \emph{point} in a topos \E is a geometric morphism
\[ p: \Set \to \E
\]
(we write $p \in \E$ by abuse of notation).
A topos \E \emph{has enough points} if given
$f \neq g: P \to Q \in \E$ there is a point $p \in \E$ such that
$p^*f \neq p^*g \in \Set$.
\end{definition}

\subsection{Stalks and germs}\label{subsec-stalks}

\begin{definition} \label{defnStalk}
\label{def:stalk}Let $(\C,J)$ be a site, $a:\, \pre{\C}\rightarrow \Sh(\C,J)$ 
the associated sheaf functor and $x \in \Sh(\C,J)$ a point. The 
\emph{stalk functor} at $x$ is given by
\[
\stal {x} \deq x^* \circ a:\, \pre {\C}\rightarrow \set.
\]
\end{definition}

Given a presheaf $F$, we say that $\stal{x}(F)$ is the stalk of $F$ at $x$. 
As an example, consider again the case $\C = \os(X)$ with
$X$ a (this time) Hausdorff topological space and $\os(X)$ its locale of opens
equipped with the open-cover topology. Let $\Sh(X)$ be the corresponding
topos of sheaves.
It can be shown
that any geometric point $x: \set \rightarrow \Sh(X)$ 
corresponds to
a ``physical'' point $x' \in X$. The stalk of 
$F \in \set^{\os(X)^{\text{op}}}$  at $x$ is then given
by 
\[
\stal{x}(F) := \colim_{U \in \os(X), x' \in U} F(U).
\]
Write $\germ_{x,U}: \; F(U) \rightarrow \stal{x}(F)$ for the canonical
map at $U$ ($\germ_x$ when $U$ is clear from the context). We call 
the equivalence class $\germ_{x,U}(s)$ of $s$ in  
$\stal{x}(F)$
the \emph{germ of $s$ at $x$}.
Obviously,
\[
\stal{x}(F) = \{ \germ_{x,U}(s) \ | \ U \in \os(X), \ x' \in U, \ s \in
F(U)\}.
\]

\subsection{Simplicial Presheaves}\label{subsec-spresheaves}

For the rest of this section, let \C be a small category with a 
Grothendieck topology $J$ such that $\Sh(\C,J)$ has \emph{enough points}.

Let $\Delta$ be the simplicial category which has objects $[n] =
\{0,1,\ldots,n\}$ for $n\geq0$, and whose morphisms are the maps such
that $x \leq y$ implies that $f(x) \leq f(y)$.
Then $\sSet$ is the category $\Set^{\Delta^{\op}}$.
This category has a well-known model structure (see~\cite{hoveyBook}
for example) where $\mathcal{W}_{\sSet}$ are the morphisms whose
geometric realization is a weak homotopy equivalence and
$\mathcal{C}_{\sSet}$ are the monomorphisms.

Objects of $\sSetCop$ are called \emph{simplicial presheaves} 
on \C since
\[
\sSetCop = \left (\Set^{\Delta^{op}} \right )^{\C^{op}}
\cong \Set^{\Delta^{op} \times \C^{op}} 
\cong \left (\Set^{\C^{op}} \right )^{\Delta^{op}}.
\]
There is an embedding\[
\begin{array}{rlcl}
 \kappa : & \pre {\mathbb{C}} & \rightarrow  & \spre {\mathbb{C}}\\
  & F & \mapsto  & \kappa _{F}\end{array}
\]
where $\kappa _{F}$ is constant \emph{levelwise} i.e. $\left(\kappa _{F}\right)\left(C\right)_{n}\deq F\left(C\right)$
for all $n\in \mathbb{N}$, and 
all the face and degeneracy maps are the identity. There is
a further embedding
\[
\begin{array}{rlcl}
 \gamma : & \sset  & \rightarrow  & \spre {\mathbb{C}}\\
  & K & \mapsto  & \gamma _{K}\end{array}
\]
where $\gamma _{K}$ is constant \emph{objectwise} i.e. $\gamma _{K}\left(C\right)\deq K$
for all $C\in \mathbb{C}$. 

Recall that for $C \in \mathbb{C}$ and $F \in \SetCop$, the Yoneda
lemma gives the isomorphism $\SetCop(y(C), F) \isom F(C)$, where $y$
is the Yoneda embedding (see Section~\ref{subsec-presheaves}).
In the simplicial case we have the following variation, which can be proved
using the same idea used in the proof of the Yoneda lemma.

\begin{proposition}
\label{pro:(Bi-Yoneda)}(Bi-Yoneda) Let $C\in \mathbb{C}$ and $F\in \sSetCop$.
There is an isomorphism
\[
\spre {\mathbb{C}}\left(\kappa _{y\left(C\right)}\times \gamma _{\Delta \left[n\right]},F\right)\cong F\left(C\right)_{n}\]
natural in all variables.
\end{proposition}

\begin{definition} \label{defnSimplicialYoneda}
Using the Yoneda embedding $y:\C \to \SetCop$ for presheaves one can
define an embedding
\[
\bar{y}: \C \xto{y} \SetCop \xto{\kappa} \sSetCop
\]
for simplicial presheaves. 
The functor $\bar{y}$ is also called a Yoneda embedding.
\end{definition}

There are two Quillen equivalent model structures on $\spre {\mathbb{C}}$
which are in a certain sense \emph{objectwise:}

\begin{itemize}
\item the \emph{projective} model structure $\spre {\mathbb{C}}_{prj}$
where $\mathcal{W}_{prj}$ and $\mathcal{F}_{prj}$ are objectwise 
(that is, $f:P \to Q \in \mathcal{W}_{prj} (\mathcal{F}_{prj})$ if and
only if for all $C \in \C$, $f(C):P(C) \to Q(C) \in
\mathcal{W}_{\sSet} (\mathcal{F}_{\sSet})$ ), 
and
\item the \emph{injective} model structure $\spre {\mathbb{C}}_{inj}$
where $\mathcal{W}_{inj}$ and $\mathcal{C}_{inj}$ are objectwise.
\end{itemize}

These were studied by Bousfield and Kan~\cite{bousfieldKanBook} and
Joyal~\cite{joyal:letter}, respectively.

\begin{proposition}
\label{pro:Both}Both $\spre {\mathbb{C}}_{prj}$ and $\spre {\mathbb{C}}_{inj}$
are proper, simplicial, cellular model categories. All objects are cofibrant
in the latter. The identity functor is a left Quillen equivalence from
the projective model structure to the injective model structure.
\end{proposition}

The injective
one is more handy when it comes down to calculating homotopical 
localizations, yet the fibrant objects are easier to grasp
in the projective one\footnote{They are objectwise Kan.}.

Using the stalk functor for presheaves, one can define a simplicial
stalk functor for simplicial presheaves.

\begin{definition} \label{defnStalkSPreSheaves}
The \emph{simplicial stalk functor} at a point $p$ in
\ShC is given by
\[
\begin{array}{rlcl}
(\_)_p: & \sSetCop & \rightarrow & \sSet
\\
& P & \mapsto & \{stalk_p(P_n)\}_{n\geq 0}.
\end{array}
\]
A morphism $f: P \to Q \in \sSetCop$ is a \emph{stalkwise equivalence} if
$f_p: P_p \to Q_p \in \sSet$ is a weak equivalence for all points $p$ in
$\ShC$.
\end{definition}

Jardine~\cite{jardineSimplicialPresheaves} proved the existence of a
local version of Joyal's injective model structure.
Since we will only be interested in the special case where \ShC has
enough points, we will not recall the definition of local weak equivalences.

\begin{theorem}[\cite{jardineSimplicialPresheaves,jardineBooleanLocalization}]
  \label{thmJardine} 
Let \C be a small category with a Grothendieck topology.
Then \sPreC the category of simplicial presheaves on \C has a proper,
simplicial, cellular model structure in which 
\begin{itemize}
\item the cofibrations are the monomorphisms, i.e. the levelwise
  monomorphisms of presheaves,
\item the weak equivalences are the \emph{local weak equivalences},
  and 
\item the fibrations are the morphisms which have the right lifting property
with respect to all trivial cofibrations.
\end{itemize}
Furthermore, if the Grothendieck topos \ShC has enough
  points, then the local weak equivalences are the stalkwise
  equivalences. 
\end{theorem}

Jardine's model structure can be seen to be cellular since it can also
be constructed as a left Bousfield localization of the injective model
structure~\cite{duggerHollanderIsaksenHypercoversaSP}.

\section{Local po-spaces} \label{sectionSummary} 

The focus of this section is to provide the reader with the main
definitions and constructions.
We define a small category of local po-spaces $\LPS$ and state some of
the properties, most of which are proved in the later sections. 
We show that Theorem~\ref{thmMain} follows from these properties and a
theorem of Jardine.


To simplify the analysis, we will only work with topological spaces
which are subspaces of $\R^n$ for some $n$, since this provides more than
enough generality for studying concurrent systems.
The main technical advantage of this setting is that we obtain small
categories.

\begin{definition}
\begin{enumerateroman}
\item 
Let \Top be the category whose objects are  subspaces
of $\R^n$ for some $n$, and whose morphisms are continuous maps.
\item
Let \Pospc be the category whose objects are \emph{po-spaces}: that is
$U \in  \Top$ together with a \emph{partial order}
(a reflexive, transitive, anti-symmetric relation) $\leq$
such that $\leq$ is a closed subset of $U \times U$ in the
product topology.
\item 
For any $M \in  \Top$ define an \emph{order-atlas} on $M$ to be an
open cover\footnote{That is, for all $i$, $U_i$ is open as a subspace
  of $M$ and $M = \cup_i U_i$.} 
$U = \{U_i\}$ of $M$ indexed by a set $I$,
where $U_i \in  \Pospc$.
These partial orders are compatible: $\leq_i$ agrees with $\leq_j$ on
$U_i \cap U_j$ for all $i,j \in I$. 
We will usually omit the index set from the notation.
\item
Let $U$ and $U'$ be two order atlases on $M$. 
Say that $U'$ is a \emph{refinement} of $U$ if
for all $U_i \in U$, and for all $x \in U_i$, there exists a $U'_j \in
U'$ such that $x \in U'_j \subseteq U_i$ and for all $a,b \in U'_j$,
$a \leq_{j'} b$ \Iff $a \leq_i b$. 
\item
Say that two order atlases are \emph{equivalent} if they have a common
refinement.
This is an equivalence relation:
reflexivity and symmetry follow from the definition.
For transitivity, if $U$ and $U'$ have a refinement $V = \{V_i\}$ and
$U'$ and $U''$ have a refinement $W = \{W_j\}$, 
let $T =  \{ V_i \cap W_j \}$.
One can check that $T$ is an order atlas of $M$ and that is a
refinement of $U'$ and $U''$.
\end{enumerateroman}
\end{definition}

Any po-space $(U,\leq)$ is a local po-space with the equivalence class
of order atlases generated by the order atlas $\{U\}$.
As a further example, we remark that any discrete space has a unique
equivalence class of order-atlases.

\begin{definition}
Let \LPS be the category of local po-spaces described as follows.
The objects, called \emph{local po-spaces}, are all pairs $(M,\coatl{U})$ 
where $M$ is an object in \Top and $\coatl{U}$ is an equivalence class
of \emph{order-atlases} of $M$. 
The morphisms, called \emph{dimaps} are described as
follows.
$f \in \LPS( (M,\coatl{U}), (N,\coatl{V}))$ \Iff $f \in \Top(M,N)$ and
for all $V = \{V_j\}_{j \in J} \in \coatl{V}$ there is a $U = \{U_i\}_{i
\in I}\in \coatl{U}$ such that for all $i \in I$, $j \in J$, for all
$x,y \in U_i \cap f^{-1}(V_j)$,
\[
x \leq_{U_i} y \implies f(x) \leq_{V_j} f(y).
\]
\end{definition}

\begin{remark}
This condition is not necessarily true for arbitrary $U \in \coatl{U}$.
For example, take $M = \{-1,1\}$ with $\coatl{U}$ the unique equivalence
class of order atlases generated by the order atlas $U = \{ \{-1\},
\{1\} \}$. Let $f = \Id_M: (M,\coatl{U}) \to (M,\coatl{U})$.
Now let $M_+$ be the po-space on $M$ with the ordering $-1 \leq 1$ and
let $M_-$ be the po-space on $M$ with the ordering $1 \leq -1$.
Then $\{M_+\} \in \coatl{U}$ and $\{M_-\} \in \coatl{U}$ (both have $U$ as
a common refinement).
However, even though $-1,1 \in M_+ \cap f^{-1}(M_-)$,
\[ 
-1 \leq_{M_+} 1 \text{ but } f(-1) \nleq_{M_-} f(1).
\]
\end{remark}

\begin{remark}
It is easy to check that a dimap of po-spaces is also a dimap of local
po-spaces.
Thus \Pospc the category of po-spaces is a subcategory of $\LPS$.
\end{remark}

\begin{remark} \label{remSubspacesInLPS}
Subobjects in \LPS. 

If $(M,\coatl{U}) \in  \LPS$, then a subspace $L \subseteq M \in \Top$
inherits local po-space structure as follows.
Let $U = \{U_i\} \in \coatl{U}$ and let $W = \{W_i\}$ where $W_i = L
\cap U_i$ and $W_i$ has the partial order inherited from $U_i$.
Then $W$ is an open cover of $L$ and the partial orders are
compatible. 
That is $W$ is an order atlas.
Let $\coatl{W}$ be the equivalence class of $W$.

We claim that $\coatl{W}$ does not depend on the choice of $U$.
Let $\tilde{U} = \{ \tilde{U}_i \} \in \coatl{U}$, let $\tilde{W}_i =
L \cap \tilde{U}_i$, and let
$\tilde{W} = \{ \tilde{W}_i \}$.
$U$ and $\tilde{U}$ have a common refinement $\hat{U} =
\{\hat{U}_i\}$.
Let $\hat{W}_i = L \cap \hat{U}_i$ and let $\hat{W} = \{\hat{W}_i\}$.
Then one can check that $\hat{W}$ is a common refinement of $W$ and
$\tilde{W}$.
So the equivalence class of $\tilde{W}$ is also $\coatl{W}$.

Next we claim that there is a dimap $\iota: (L, \coatl{W}) \to
(M,\coatl{U})$ given by the inclusion $\iota: L \incl M$.
Let $U = \{U_k\} \in \coatl{U}$, let $W_k = L \cap U_k$, and let
$W = \{W_k\}$. Then $W \in \coatl{W}$.
Let $x,y \in W_j \cap \iota^{-1}(U_k) = W_j \cap L \cap U_k
= W_j \cap W_k.$
Note that $\iota(x) = x$ and $\iota(y) = y$.
Then
\[
x \leq_{W_j} y \iff x \leq_{W_k} y \iff x \leq_{U_k} y.
\]
Therefore when $L \subseteq M \in \Top$, then there is an induced
inclusion $(L,\coatl{W}) \subseteq (M, \coatl{U}) \in \LPS$.
\end{remark}

The remark above will be used implicitly and without reference in
Section~\ref{sectionPoints}.

\begin{definition}
A collection of dimaps $\{ \phi_j: (M_j, \coatl{U}^j) \to
(M,\coatl{U}) \}$ $\LPS$ is an \emph{open dicover} if 
\begin{enumerateroman}
\item
$\{\phi_j: M_j \to M\}$ is an open cover, and 
\item
for each $j$, 
$\coatl{U}^j$ is the local po-space structure inherited from $(M,\coatl{U})$.
\end{enumerateroman}
\end{definition}

\begin{remark} \label{remSubPoSpaces}
The local po-space structures inherited by the subspaces of
$(M,\coatl{U})$ are compatible.
So if $\{\phi_j: (M_j,\coatl{U}^j) \to (M,\coatl{U})\}$ is a open cover,
then for each $j$, there is a $U^j = \{U^j_k\} \in \coatl{U}^j$ such that
$U' = \{U^j_k\}_{j,k}$ is an order atlas for $M$ and $U' \in \coatl{U}$.
\end{remark}

The following 
is easy to check.

\begin{lemma} \label{lemmaSmall}
\Top and \LPS are small categories.
\end{lemma}

Define $U: \LPS \to \Top$ to be the forgetful functor defined on
objects and morphisms as follows $(M,U) \mapsto M$ and $\varphi \mapsto
\varphi$.

Define $F: \Top \to \LPS$ as follows. 
If $M$ is an object in \Top, then let $F(M) =
(M,\Rmphi)$, where $\Rmphi$ is the equivalence class of $M_{\phi} = \{
M \}$ with $x \leq_M y \iff x=y$.
If $f: M \to N \in \Top$,
then $F(f) = f: (M, \Rmphi) \to (N, \Rnphi)$.
This is a dimap since for any $V = \{V_j\} \in \Rnphi$ with $x,y
\in f^{-1}V_j$, $x \leq_M y \implies x=y \implies
f(x)=f(y) \implies f(x) \leq_{V_j} f(y)$.

\begin{remark}
Note that $U$ is faithful and $F$ includes \Top as a full subcategory
of \LPS.
\end{remark}

\begin{proposition}
$F: \Top \rightleftarrows \LPS: U$ is an adjunction.
\end{proposition}

\begin{proof}
Let $M$ be an object in $\Top$ and $(N,\bar{V}) \in  \LPS$.
We claim that there is a natural bijection
\[
\LPS( F(M), (N,\bar{V}) ) \isom \Top (M, U(N,\bar{V})).
\]
We need to show that there is a natural bijection
\begin{equation*} \label{eqnTopLPSadjunction}
\theta: \Top(M,N) \isomto \LPS( (M, \Rmphi), (N,\bar{V}) ).
\end{equation*}

If $f \in \LPS((M,\Rmphi), (N,\bar{V}))$, then $f \in \Top(M,N)$ such
that for any $V = \{V_j\} \in \bar{V}$, for
all $j$, $f|_{f^{-1}({V}_j)}$ satisfies
$x \leq_M y \implies f(x) \leq_{{V}_j} f(y)$.
Since $x \leq_M y$ if and only if $x=y$ this last condition is vacuous.
Thus the bijection is simply $\theta: f \mapsto f$.

To show naturality let $\alpha: (N,\bar{V}) \to (N',\bar{V'}) \in \LPS$
and $\xi: M' \to M \in \Top$.
Then 
\[ \theta (U (\alpha) \circ f \circ \xi) = \alpha \circ f \circ \xi =
\alpha \circ \theta(f) \circ \xi.
\]
\end{proof}

\begin{remark}
\LPS does not have colimits.

Consider the product of the directed circle and an interval.
Now collapse the top circle of this cylinder.
The vertex of the resulting cone does not have a local partial order.
\end{remark}

\section{The open-dicover topology}
\label{sectionSubcanonical} 

We define the open cover Grothendieck topology for \Top and the open
dicover Grothendieck topology for \LPS in the following lemma.
The proof of the lemma follows directly from the definition of a basis
for a Grothendieck topology.

\begin{lemma} \label{lemmaOpenCoverTopology}
\begin{enumerate}
\item
\Top has a \emph{Grothendieck topology} whose basis is given by
the open covers.
For $M \in  \Top$ let $K(M) = \{\text{open covers of }M\}$.
Let $J$ be the Grothendieck topology generated by $K$. 
Call $J$ the \emph{open cover topology}. 
\item
Analogously, \LPS has a Grothendieck topology whose basis is given by
the {open dicovers in $\LPS$}.
Let $K((M,\coatl{U})) = \{\text{open dicovers of } (M,\coatl{U})\}$.
Call the Grothendieck topology generated by $K$ the \emph{open-dicover
topology}.
\end{enumerate}
\end{lemma}

In Section~\ref{sectionSummary}, we defined a Grothendieck topology to
be subcanonical if every representable presheaf a sheaf.
In this section, we will prove that the open-dicover topology is
subcanonical.

The following proposition shows that if a Grothendieck topology is
generated by a basis $K$, then to see if a presheaf is a sheaf it
suffices to check the basis. For the definition of matching families
and amalgamations see Remark~\ref{defnMatchingFamily}.

\begin{proposition}[{\cite[Proposition III.4.1]{\mmBook}}]
  \label{propPresheafIsSheaf}
Let \C be a small category with a Grothendieck topology $J$ generated
by a basis $K$.
Then a presheaf $P \in \SetCop$ is a sheaf for $J$ if and only if for
every $M \in  \C$ and every cover $\{\phi_j: M_j \to M\}
\in K(M)$, every matching family for $\{\phi_j\}$ of elements of $P$
has a unique amalgamation. 
\end{proposition}

\begin{example} \label{egOpenCoverTopologySubcanonical}
Let $N \in  \Top$ and $y(N) = \Top(-,N) \in  \Set^{\Top^{\op}}$.
Let $\phi_j: M_j \to M$ be an open cover, and let $\alpha_j: M_j \to
 N$ be a matching family.
Then $\phi_j$ has a unique amalgamation $\phi: M \to N$.
Therefore $y(N)$ is a sheaf for the open cover topology, and hence the
 open cover topology is subcanonical. 
\end{example}

\begin{proposition} \label{propDicoverTopSubcanonical}
 In the open-dicover topology $J$ for local po-spaces every
 representable presheaf is a sheaf. That is $J$ is subcanonical.
\end{proposition}

\begin{proof}
Consider the representable presheaf 
\[
y((N,\bar{V})) = \LPS(-,(N,\bar{V})) \in \SetLPSop.
\] 
By Proposition~\ref{propPresheafIsSheaf}, $y((N,\bar{V}))$ is a sheaf
\Iff for all open dicovers $\{ \phi_j \} \in K( (M,\bar{U}) )$, 
any matching family 
\[
\{\alpha_j: (M_j,\bar{U}_j) \to (N,\bar{V})\}
\]
has a unique amalgamation $\alpha: (M,\bar{U}) \to (N,\bar{V})$.
That is, there is a map $\alpha$ such that the diagrams
\[
\xymatrix{
(M_j, \bar{U}_j) \ar[r]^{\phi_j} \ar[d]_{\alpha_j} & (M,\bar{U})
\ar[dl]^{\alpha} \\
(N,\bar{V})
}
\]
commute in \LPS for all $j$.

Let $\{\alpha_j\}$ be such a matching family for an open dicover
$\{\phi_j\}$.
Since $\{\phi_j\}$ is an open dicover, then by
Remark~\ref{remSubPoSpaces} for each $j$ there is a $U^j = 
\{U^j_k\} \in \bar{U}_j$ such that $U' = \{U^j_k\}_{j,k}$ is an order
atlas and $U' \in \bar{U}$.

By definition $\{\phi_j: M_j \to M\}$ is a cover in \Top and
$\{\alpha_j: M_j \to N\}$ is a matching family.
Therefore there is a unique amalgamation $\alpha: M \to N \in \Top$.
That is, there is a map $\alpha$ such that
\[
\xymatrix{
M_j \ar[r]^{\phi_j} \ar[d]_{\alpha_j} & M \ar[dl]^{\alpha} \\
N
}
\]
commutes in \Top for all $j$.
It remains to show that $\alpha$ is a dimap.
Let $V = \{V_l\} \in \bar{V}$.
Since $\alpha_j: (M_j, \bar{U}_j) \to (N,\bar{V}) \in \LPS$, there is
a $\tilde{U}^j = \{\tilde{U}^j_k\}_k \in \bar{U}^j$ such that for all $k,l$,
\[
\text{for all }x,y \in \tilde{U}^j_k \cap \alpha^{-1}_j(V_l), \quad x
\leq_{\tilde{U}^j_k} y \implies \alpha_j(x) \leq_{V_l} \alpha_j(y).
\]
Now for each $j$, let $\hat{U}^j = \{\hat{U}^j_k\}_k \in \bar{U}_j$ be
a common refinement of $\tilde{U}^j$ and $U^j$.
Then since $\hat{U}^j$ is a refinement of $\tilde{U}^j$,
\begin{equation} \label{eqnAlphaJMonotone}
\text{for all }x,y \in \hat{U}^j_k \cap \alpha^{-1}_j(V_l), \quad x
\leq_{\hat{U}^j_k} y \implies \alpha_j(x) \leq_{V_l} \alpha_j(y),
\end{equation}
and since $\hat{U}^j$ is a refinement of $U^j$, if we define $U = \{
\hat{U}^j_k \}_{j,k}$, then $U \in \bar{U}$.

Since $\alpha$ is an amalgamation of $\{\alpha_j\}$ in \Top
if $x \in \hat{U}^j_k \subset M$, then $\alpha(x) = \alpha_j(x)$ and
for all $l$,
$\hat{U}^j_k \cap \alpha_j^{-1}(V_l) = \hat{U}^j_k \cap
\alpha^{-1}(V_l)$.
Therefore using~\eqref{eqnAlphaJMonotone} for all $k,l$,
\begin{equation*} \label{eqnAlphaMonotone}
\text{for all }x,y \in \hat{U}^j_k \cap \alpha^{-1}(V_l), \quad x
\leq_{\hat{U}^j_k} y \implies \alpha(x) \leq_{V_l} \alpha(y).
\end{equation*}
That is $\alpha$ is a dimap.
Therefore $\alpha: (M,\bar{U}) \to (N,\bar{V})$ is a unique
amalgamation of $\{\alpha_j\}$.
\end{proof}

\section{Equivalence of sheaves and di-\'etale bundles}

In this section \C is either \Top or \LPS with the Grothendieck topology
generated by open (di)covers.

\begin{notation}
We will use $A \opensubset B$ to denote that $A$ is an open subset of
$B$.
\end{notation}

\begin{notation}
Let $Z \in  \C$ and let $F \in \SetCop$. Choose $x \in U \opensubset
Z$ and $s \in F(U)$.
Then for open subobjects of $U$, $L \stackrel{i}{\incl} U$, we have
$F(i): F(U) \to F(L)$ 
and we will use the notation
\[
s|_L := F(i)(s).
\]
Recall that $\stal{x}(F) = \colim_{x \in L \opensubset U} F(L)$ and
$\germ_x(s)$ is the equivalence class represented by $s$ in $\stal{x}(F)$.
\end{notation}

\begin{definition}
Given $Z \in  \C$, a \emph{bundle} over $Z$ is just a morphism $p:W
\to Z \in \C$. An \emph{(di)\'{e}tale bundle} is a bundle which
is a \emph{local (di)homeomorphism}. 
That is, given $y \in W$ there is some open set $V \subset W$ such that
$p(V)$ is open in $Z$ and $p|_V$ is an isomorphism in $\C$.

A morphism of (\'{e}tale) bundles $p:W \to Z$ and $p:W' \to Z$ is a
morphism $\theta: W \to W' \in \C$ such that the following diagram commutes:
\[
\xymatrix{
W \ar[rr]^{\theta} \ar[dr]_{p} & & W' \ar[dl]^{p'} \\
& Z
}
\]

Let $\EtaleZ$ denote the category of (di)\'{e}tale bundles over $Z$.
In addition let $\OZ$ denote the category of open subobjects of
$Z$, where the objects are open subobjects of $Z$ and the morphisms are
the inclusions.
\end{definition}

\begin{theorem}[Theorem \ref{thmEtaleSheafCorrespondence}]
Let $Z \in  \C$. Then there is an equivalence of categories
\[
\Gamma: \EtaleZ \adjn \ShOZ: \Lambda.
\]
\end{theorem}

\begin{proof}
It is well known that the statement of
Theorem~\ref{thmEtaleSheafCorrespondence} is true when $\C = \Top$
(see for example~\cite[Corollary II.6.3]{\mmBook}).
We will show that this equivalence between {\'e}tale bundles on
topological spaces and sheaves on topological spaces extends to local
po-spaces.

First we describe the functors $\Gamma$ and $\Lambda$ in the case
where $\C = \Top$.
The functor $\Gamma$ assigns to each bundle $W \xto{p} Z$ the
presheaf of cross-sections: 
\begin{eqnarray*}
P: \OZ^{\op} & \to & \Set \\
U & \mapsto & \{ s: U \to W \in \C \ | \ p \circ s = \Id_U \} \\
U \stackrel{\theta}{\incl} V & \mapsto & \theta^* \quad (\theta^*(t) =
t \circ \theta).
\end{eqnarray*}
One can check that if $p$ is {\'e}tale, then $P$ is in fact a
sheaf~\cite[p.79]{\mmBook}. 
Thus $\Gamma$ restricts to a functor $\Gamma: \EtaleZ \to \ShOZ$.

Given a presheaf $P: \OZ^{\op} \to \Set$, $\Lambda(P)$ is the bundle
$W \xto{p} Z$ where
\[
W = \{ \germ_x s \ | \ x \in U \opensubset Z, s \in P(U)\} \text{ and }
p: \germ_x s \mapsto x.
\]
A basis for the topology on $W$ is given by the sets $\dot{s}(U)$,
where $U$ is an open set in $Z$, $s \in P(U)$ and 
\begin{eqnarray*}
\dot{s}: U & \to & \Lambda(P) \\
x & \mapsto & \germ_x s.
\end{eqnarray*}
Using this topology, $p:W \to Z$ is a continuous map.
Again, one can check that if $P$ is a sheaf, then $W \xto{p} Z$ is in
fact an {\'e}tale bundle~\cite[p.85]{\mmBook}. 
So $\Lambda$ restricts to a functor $\Lambda: \ShOZ \to \EtaleZ$.

Now we will show that $\Gamma$ and $\Lambda$ can be similarly defined
in the case where $\C = \LPS$.
Let $p: (W,\bar{T}) \to (Z,\bar{U})$ be an {\'e}tale bundle of local
po-spaces. 
The definition of $\Gamma$ is exactly the same: $\Gamma((W,\bar{T})
\xto{p} (Z,\bar{U}))$ is the sheaf of cross-sections.

Given a sheaf $P$ on a local po-space $(Z,\bar{U})$, $\Lambda(P) = (W
\xto{p} Z)$ is an {\'e}tale bundle of topological spaces.
To extend $\Lambda$ to local po-spaces it remains to define a local
order on $W$ and show that this makes $p$ a dimap.

\begin{lemma} \label{lemmaCanonicalLPSstructure}
$W$ has a canonical local po-space structure such that $p$ is a dimap.
\end{lemma}

\begin{proof}

Recall that the sets $\dot{s}(U)$ defined above form a basis for the
topology of $W$. 
Choose an order atlas $\{(U_i,\leq_i)\} \in \bar{U}$ for $Z$.
For each open sub-po-space $V \subset U_i$ and each $s \in P(V)$, 
$\dot{s}(V) \subset W$ is a po-space under the relation 
\[ \germ_x s \leq_{\dot{s}(V)} germ_y s \text{ if and only } x \leq_i y.
\]
This is well-defined since $\{U_i\}$ is an order-atlas, and it makes
$\dot{s}(V)$ a po-space since $\dot{s}: U_i \to \dot{s}(U_i)$ is a
homeomorphism. 

We claim that 
\[ T := \{\dot{s}(V) \ | \ V \opensubset U_i, s \in P(V)\}
\]
is an order atlas on $W$.
First we need to show that it is an open cover.
Each of the sets is open by construction.
If $U \in \OZ$ and $s \in P(U)$, consider $\germ_x s$.
Since $\{U_i\}$ is an open cover of $Z$, for some $i$, $x \in U_i$.
Let $V = U \cap U_i$.
Then $\germ_x s = \germ_x s|_V \in \dot{(s|_V)}(V)$.
Therefore $T$ is an open cover of $W$.

Finally we need to show that the orders are compatible.
For $k=1,2$ let $V_k \opensubset U_{i_k} \opensubset Z$, and $s_k \in P(V_k)$.
Assume $g_1, g_2 \in \dot{s_1}(V_1) \cap \dot{s}(V_2)$.
That is, $g_1 = \germ_{x_1} s_1$ $= \germ_{x_1} s_2$ and 
$g_2 = \germ_{x_2} s_1 = \germ_{x_2} s_2$.
For $k=1,2$, 
\[
g_1 \leq_{\dot{s_k}(V_k)} g_2 \iff x_1 \leq_{i_k} x_2.
\]
Since $\{U_i\}$ is an order-atlas, the order $\leq_{i_1}$ and
$\leq_{i_2}$ are compatible.
Therefore the orders $\leq_{\dot{s_1}(V_1)}$ and
$\leq_{\dot{s_1}(V_1)}$ are compatible, and $T$ is an
order-atlas on $W$.

Let $\bar{T}$ be the equivalence class of order atlases of $T$.
We claim that $\bar{T}$ does not depend on the choice of $U \in
\bar{U}$.

Let $U, U' \in \bar{U}$, then $U$ and $U'$ have a common refinement
$U''$.
Let $T, T', T''$ be the corresponding order-atlases for $W$
constructed as above. 
We will show that $T''$ is a refinement of $T$.

Let $A \opensubset U_j \in U$, $s \in P(A)$ and $\germ_x s \in
\dot{s}(A)$.
Then there is some $U''_k \in U''$ such that $x \in U''_k$ and $U''_k$
is a sub-po-space of $U_j$.
Let $A'' = A \cap U''_k$.
It follows that $\dot{(s|_{A''})}(A'') \subset \dot{s}(A)$, and 
$\germ_x s = \germ_x (s|_{A''}) \in \dot{(s|_{A''})}(A'') \in T''$.
Since $U''_k$ is a sub-po-space of $U_j$ it follows that 
$\dot{(s|_{A''})}(A'')$ is a sub-po-space of $\dot{s}(A)$.
Thus $T''$ is a refinement of $T$.

Similarly $T''$ is a refinement of $T'$ and is hence a common
refinement of $T$ and $T'$. 
Therefore $\bar{T}$ does not depend on the choice of $U \in \bar{U}$.

Finally we will show that the projection $p:W \to Z$ given by $\germ_x s
\mapsto x$ is a dimap.
Let $U \in \bar{U}$ be an order-atlas on $Z$. 
Let $T$ be the order-atlas on $W$ constructed above from $U$.
Observe that $T \in \bar{T}$, since $\bar{T}$ does not depend on the
choice of $U \in \bar{U}$.
Let $U_j \in U$, let $A \opensubset U_i \in U$, and let $s \in P(A)$.
Assume that
\[ \germ_{x_1} s, \germ_{x_2} s \in \dot{s}(A) \cap p^{-1}(U_j).
\]
Then $x_1, x_2 \in U_i \cap U_j$.
By the construction of $T$ and since $U$ is an order atlas,
\[
\germ_{x_1} s \leq_{\dot{s}(A)} \germ_{x_2} s \iff x_1 \leq_{U_i} x_2
\iff x_1 \leq_{U_j} x_2.
\]
Therefore $\Lambda$ can be extended to local po-spaces.
\end{proof}

Thus we have maps
\[
\Gamma: \EtaleZ \adjn \ShOZ: \Lambda.
\]
To show that they give an equivalence of categories we will show that
for a sheaf $P$ and an {\'e}tale space $W \xto{p} Z$ there are natural
isomorphisms 
\[
\epsilon_W: \Lambda \Gamma W \to W \text{ and } \eta_P: P \to \Gamma
\Lambda P.
\]
Recall that elements of $\Lambda \Gamma W$ are of the form $\dot{s}(x)
= \germ_x s$, where $s:U \to W$ satisfies $p \circ s = \Id_U$ and $x \in U$.
Define $\epsilon_W$ to be the map $\dot{s}x \mapsto sx$.
We will show this is an isomorphism by constructing an inverse
$\theta_W$.
Let $y \in W$ and let $x = py$.
Since $W$ is {\'e}tale there exists $y \in V \opensubset W$ such that
$p|_V: V \isomto p(V)$.
Let $q = (p|_V)^{-1}$.
Then define $\theta_W(y) = \germ_x q = \dot{q}x$.
Then we claim $\theta_W$ is an inverse for $\epsilon_W$.
Indeed
\[
\epsilon_W \theta_W y = \epsilon_W \dot{q}x = qx = y.
\]
Also for all $\dot{s}x \in \Lambda \Gamma W$,
$\theta_W \epsilon_W \dot{s}x = \theta_W sx = \germ_x t$,
where t is a restriction of $s$.
So $\germ_x t = \germ_x s = \dot{s}x$.

Finally we claim that $\epsilon_W$ and $\theta_W$ are dimaps.
First choose $T = \{T_k\} \in \bar{T}$ and $U = \{U_i\} \in \bar{U}$
such that $p$ satisfies the dimap condition.
>From $T$ construct the canonical order atlas of the form
$\{\dot{s}V\}$ for $\Lambda \Gamma W$ as in the proof of
Lemma~\ref{lemmaCanonicalLPSstructure}. 
Now let $\dot{s}x_1, \dot{s}x_2 \in \dot{s}V \cap
\epsilon_W^{-1}(T_k)$.
Then by construction,
\[
\dot{s}x_1 \leq_{\dot{s}V} \dot{s}x_2 \iff x_1 \leq_{U_i} x_2.
\]
Since $s$ satisfies the dimap condition this implies that
$sx_1 \leq_{T_k} sx_2$ which is the same as $\epsilon_W \dot{s}x_1
\leq_{T_k} \epsilon_W \dot{s}x_2$.
Thus $\epsilon_W$ is a dimap.
Next let $y_1, y_2 \in T_k \cap \theta_W^{-1}(\dot{s}V) = T_k \cap
\epsilon_W(\dot{s}V) = T_k \cap sV$.
Then there are $x_1,x_2 \in V$ such that $y_1 = sx_1$ and $y_2 =
sx_2$.
Since $p$ satisfies the dimap condition
\[
y_1 \leq_{T_k} y_2 \implies py_1 \leq_{U_i} py_2.
\]
But this is the same as $x_1 \leq_{U_i} x_2$ which implies that
$\dot{s}x_1 \leq_{\dot{s}V} \dot{s}x_2$.
Therefore $\theta_W$ is a dimap.

The proof that the morphism $\eta_P$ is a bijection is the same as the
proof in the case of topological spaces~\cite[Theorem II.5.1]{\mmBook}.
\end{proof}

\section{Points} \label{sectionPoints}

In this section \C is either \Top or \LPS with the Grothendieck topology
generated by open (di)covers.

Let \SetCop and \ShC be the topoi of presheaves and sheaves on $\C$.
Recall that the inclusion functor $i: \ShC \to \SetCop$ has a right
adjoint $a$ called the associated sheaf functor.
Recall from Definition~\ref{defnStalk} that 
if $p$ is a point in \ShC and $\alpha \in \SetCop$, then 
$\stal{p}(F) = p^* \circ a(\alpha)$.

Let $Z \in  \C$.
Then $Z$ is a topological space or a local po-space and we can choose
any point (in the usual sense) $x \in Z$. 
Define
\begin{align*}
p_x^*:  \Set^{\C^{\op}} & \to \Set \\
 F & \mapsto \colim_{x \in L \opensubset Z} F(L)
\end{align*}
where the colimit is taken over all open subsets of $Z$ containing $x$.
See Remark~\ref{remSubspacesInLPS} for a discussion of subobjects in
$\LPS$. 

Given a functor $p^*: \SetCop \to \Set$ there is an induced functor
\[ A: \C \xto{y} \SetCop \xto{p^*} \Set,
\]
where $y$ is the Yoneda embedding defined on objects and morphisms by 
$Z \mapsto \C(-,Z)$ and $\varphi \mapsto \C(-,\varphi)$.

Given a functor $A: \C \to \Set$ one can define induced adjoint functors
$p^*: \SetCop \!\! \to \Set$ and $p_*: \Set \to \SetCop$ 
($p^* = -\tensor_{\C} A$ and $p_* = \C(A,-)$, see \cite[Section
VII.2]{\mmBook} ).

\begin{definition}
\begin{enumerateroman}
\item
The functor $A: \C \to \Set$ is \emph{flat} if the corresponding
$p^*$ is left exact. 
\item
$A$ is \emph{continuous} if $A$ sends each covering sieve to an epimorphic
family of functions.
That is, if $S$ is a covering sieve, then the family of functions
$\{A(\varphi) | \varphi \in S\}$ is jointly surjective.
\end{enumerateroman}
\end{definition}

\begin{proposition}[{\cite[Corollary VII.5.4]{\mmBook}}]
Using the correspondence above, $p$ is a point in $\SetCop$ if and
only if $A$ is flat.
Furthermore $p$ descends to a point in $\ShC$ if and only if $A$ is flat and
continuous.
\end{proposition}

\begin{proposition}
$p_x$ defined above descends to a point in \ShC.
\begin{equation} \label{diagPx}
\xymatrix{
\SetCop \ar[r]^{p_x^*} \ar@<0.5ex>[d]^a & \Set \\
\ShC \ar@<0.5ex>[u]^i \ar@{-->}[ur]
}
\end{equation}
\end{proposition}

\begin{proof}
Let $A_x = p_x^* \circ y$, where $y$ is the Yoneda embedding.

First we show that $p_x^*$ is left exact, that is it preserves finite
limits.
Let $F \times_G H$ be a pullback in \SetCop.
\begin{eqnarray*}
p_x^* (F \times_G H) & = & \colim_{x \in L \subseteq Z} (F \times_G
H)(L) \\
& = & \colim_{x \in L \subseteq Z} F(L) \times_{G(L)} H(L) \\
& = & \colim F(L) \times_{\colim G(L)} \colim H(L) \\
& = & p_x^* F \times_{p_x^* G} p_x^* H
\end{eqnarray*}
The third equality holds because $\colim$ commutes with pullbacks in
$\Set$, and the others are by definition.
Thus $A$ is flat and $p_x^*$ is a point in \SetCop.

Next we show that $A_x$ is continuous.
Let $\{Y_i \xto{\varphi_i} N\}$ be a covering sieve for $N$ in \C.
Recall that $A_x = p_x^* \circ y$.
Let $(\varphi_i)_*$ denote composition with $\varphi_i$.
For each arrow in the covering sieve,
\begin{eqnarray*}
p_x^* \circ y(Y_i \xto{\varphi_i} N) & = & p_x^*( \C(-,Y_i)
\xto{(\varphi_i)_*} \C(-,N) ) \\
& = & \colim_{x \in L \subseteq Z} ( \C(L,Y_i) \xto{(\varphi_i)_*} \C(L,N)
) \\
& = & y(Y_i)_x \xto{(\varphi_i)_*} y(N)_x.
\end{eqnarray*}

We claim that this is an epimorphic family of functions in \Set.
Let $f \in y(N)_x$.
Then there is an open subspace $L$ such that $x \in L \subseteq Z$ and $f$ is
represented by a morphism $f' \in \C(L,N)$.
Since $\{Y_i\}$ covers $N$, $f'(x) \in Y_k$ for some $k$.
Let $K = (f')^{-1}(Y_k)$.
Then $K$ is open and $x \in K \subseteq L \subseteq Z$.
Furthermore $f'|_K \in \C(K,Y_k)$ which represents an element 
$f'' \in y(Y_k)_x$, and $(\varphi_k)_* f'' = f$.
Hence we have an epimorphic family as claimed.
Thus $A$ is continuous and $p_x$ descends to a point in \ShC.
\end{proof}

Abusing notation we will also denote the induced functor in
diagram~\eqref{diagPx} by $p_x^*$.
With this abuse of notation, the stalk of $F \in \SetCop$ at $x$
is given by $\stal{x}(F) = p_x^* a (F) = p_x^* (F)$.
Note that $\stal{x}(F) = \{\germ_x(s) \ | \ x \in U \opensubset Z, \ s
\in F(U)\}$.

\begin{theorem} \label{thmEnoughPoints}
The points $p_x$ defined above provide enough points for $\ShC$.
That is, given $f \neq g: P \to Q \in \ShC$, there is an $Z \in 
\C$ and a $x \in Z$ such that $p_x^*f \neq p_x^*g: p_x^*P \to p_x^*Q
\in \Set$.
\end{theorem}

\begin{proof}
Given $Z \in  \C$ and either $P \in  \ShC$ or $f \in \Mor \ShC$,
let $P_Z$ or $f_Z$ denote the restriction to $\ShOZ$.

Assume that $f \neq g: P \to Q \in \ShC$.
Thus there is some $Z \in  \C$ such that $f_Z \neq g_Z: P_Z \to Q_Z
\in \ShOZ$.

By Theorem~\ref{thmEtaleSheafCorrespondence}, this is equivalent to
saying that the corresponding maps between {\'e}tale spaces are not
equal.
That is,
\[
\Lambda f_Z \neq \Lambda g_Z: \Lambda P_Z \to \Lambda Q_Z \in \EtaleZ.
\]
Thus there is some point $y \in \Lambda P_Z$ such that $\Lambda f_Z(y)
\neq \Lambda g_Z (y)$.

By the definition of $\Lambda$, $y = \germ_x s$ for some $x \in
U \opensubset Z$ and $s \in P_Z(U)$.
That is $y \in \stal{x}(P) = p_x^* P$.
Therefore $p_x^* f \neq g_x^* g: p_x^* P \to p_x^* Q$. 
\end{proof}

\section{Stalkwise equivalences} \label{sectionStalkwiseC}

Let $(\C, \tau)$ be a site with a subcanonical Grothendieck
topology such that \ShC has enough points and let $\bar{y}: \C \to
\sSetCop$ be the Yoneda embedding. 
Recall the definition of stalkwise equivalence in
Definition~\ref{defnStalkSPreSheaves} which uses the simplicial stalk
functor $(\cdot)_p$.
Also recall the Yoneda embedding $\bar{y}: \C \to \sSetCop$ given in
Definition~\ref{defnSimplicialYoneda}. 
Let $\varphi: X \to Y \in \C$.

\begin{lemma} \label{lemmaSetIso}
$\bar{y} (\varphi)$ is a stalkwise equivalence if and only if for all
points $p$ in $\ShC$, $p^* a y (\varphi) \in \Set$ is an isomorphism.
\end{lemma}

\begin{proof}
Let $p$ be a point in \ShC.
Recall that the simplicial stalk of $\bar{y}(\varphi)$ at $p$ is given by
\[
(\bar{y} (\varphi))_p = \{\stal{p}(\bar{y}(\varphi)_n)\}_{n\geq 0}
= \{ p^* a y (\varphi) \}_{n\geq 0},
\]
which is simplicially constant.
Thus $\bar{y} (\varphi)_p \in \sSet$ is an isomorphism if and only if
$p^* a y (\varphi) \in \Set$ is an isomorphism.
\end{proof}

\begin{lemma} \label{lemmaAY}
If the Grothendieck topology $\tau$ is subcanonical, then the composite functor $\C \xto{y}
\SetCop \xto{a} \ShC$ is faithful.
\end{lemma}

\begin{proof}
By the Yoneda lemma, $y$ is full and faithful.
Since $\tau$ is subcanonical $\im(y) \subset \ShC$.
Furthermore $a \circ i: \ShC \to \ShC$ is
naturally isomorphic to the identity functor~\cite[Corollary
III.5.6]{\mmBook}. 
Thus $a y$ is naturally isomorphic to $y$ which is faithful.
\end{proof}

\begin{theorem} \label{thmBijection}
Let $\varphi: X \to Y \in \C$ and assume that $\bar{y} (\varphi)$ is a
stalkwise equivalence. 
Then $\varphi$ is bijective.
\end{theorem}

The proof of this theorem is split into the following two propositions.

\begin{proposition}
Let $\varphi: X \to Y \in \C$ and assume that $\bar{y} (\varphi)$ is a
stalkwise equivalence. 
Then $\varphi$ is epi.
\end{proposition}

\begin{proof}
For $i=1,2$, let $\psi_i: Y \to Z \in C$ be a morphism such that $\psi_1
\circ \varphi = \psi_2 \circ \varphi: X \to Z$.
Then for all points $p$ in \ShC, $p^* a y (\psi_1 \circ \varphi) = p^*
a y (\psi_2 \circ \varphi)$.
>From this it follows that 
\[
p^* a y (\psi_1) \circ p^* a y (\varphi) = p^*
a y (\psi_2) \circ p^* a y (\varphi).
\]
But by Lemma~\ref{lemmaSetIso} $p^* a y (\varphi)$ is a set
isomorphism, so in particular it is epi.
Therefore $p^* a y \psi_1 = p^* a y \psi_2$ for all points $p$ in
$\ShC.$
Since $\C$ has enough points, $a y \psi_1 = a y \psi_2$.
By Lemma~\ref{lemmaAY} $a \circ y$ is faithful, thus $\psi_1 = \psi_2$.
Therefore $\varphi$ is epi.
\end{proof}

\begin{proposition}
Let $\varphi: X \to Y \in \C$ and assume that $\bar{y} (\varphi)$ is a
stalkwise equivalence. 
Then $\varphi$ is mono.
\end{proposition}

\begin{proof}
For $i=1,2$, let $\psi_i: W \to X \in C$ be a morphism such that
$\varphi \circ \psi_1 = \varphi \circ \psi_2: W \to Y$.
As in the proof of the previous proposition, for all points $p$ in \ShC, 
\[
p^* a y (\varphi) \circ p^* a y (\psi_1) = p^*
a y (\varphi) \circ p^* a y (\psi_2).
\]
Again by Lemma~\ref{lemmaSetIso}, $p^* a y (\varphi)$ is mono.
Therefore $p^* a y \psi_1 = p^* a y \psi_2$ for all points $p$ in
\ShC.
Since $\C$ has enough points, $a y \psi_1 = a y \psi_2$.
By Lemma~\ref{lemmaAY} $a \circ y$ is faithful, thus $\psi_1 = \psi_2$.
Therefore $\varphi$ is mono.
\end{proof}


Let \C = \Top or \LPS with the open cover topology.
By Example~\ref{egOpenCoverTopologySubcanonical} and
Proposition~\ref{propDicoverTopSubcanonical} this topology is 
subcanonical.

Recall from Section~\ref{sectionPoints} that if $Z \in  \C$ and $x
\in Z$, then 
\begin{equation} \label{eqnPx}
\begin{split}
p_x^*:  \Set^{\C^{\op}} & \to \Set \\
 F & \mapsto \colim_{x \in L \opensubset Z} F(L)
\end{split}
\end{equation}
descends to a point in \ShC (where the colimit is taken over open
subspaces of $Z$ which contain $x$).

\begin{theorem} \label{thmStalkwiseTopLPS}
Let $\varphi: X \to Y \in \C$.
Then $\bar{y} (\varphi)$ is a stalkwise equivalence if and only if $\varphi$
is an isomorphism in $\C$. 
\end{theorem}

\begin{proof}
($\Leftarrow$)
If $\varphi$ is an isomorphism, then for all points $p$ in \ShC 
$p^*ay (\varphi)$ is an isomorphism.
Hence by Lemma~\ref{lemmaSetIso} $\bar{y} (\varphi)$ is a stalkwise
equivalence.

\noindent
($\Rightarrow$)
Assume that $\bar{y} (\varphi)$ is a stalkwise equivalence.
Then by 
Theorem~\ref{thmBijection}, $\varphi$ is a bijection.

Let $x \in Y$.
Let $p_x$ be the corresponding point defined in~\eqref{eqnPx}.
Then 
\[
p_x^*ay (\varphi): \colim_{x \in L \subseteq Y} \C(L,X) \xto{\varphi_*}
\colim_{x \in L \subseteq Y} \C(L,Y) \in \Set
\]
is a bijection.
Let $f: Y \to Y$ be given by $f = \Id_Y$.
Let $\bar{f} = [f] \in \colim_{x \in L \subseteq Y} \C(L,Y)$.
Let $\bar{g} = (\varphi_*)^{-1}(\bar{f})$.
Then there is some $x \in W \subseteq Y$ such that $\bar{g}$ has a
representative $g \in \C(W,X)$.

Let $f' = \varphi_* g = \varphi \circ g$.
Then $[f'] = \varphi_*[g] = [f]$.
Therefore there exists $x \in S \subseteq Y$ such that $S \subset Y
\cap W$ and $f'|_S = f|_S = \Id_Y|_S$.

Let $\psi = g|_S$.
Therefore $\varphi \psi = \Id_S$.
Let $T = \im(\psi)$.
Then $\varphi|_T \circ \psi = \Id_S$ and $\varphi|_T$ is a bijection.
Hence $\varphi|_T: T \to S$ is an isomorphism, where $x \in S$.

Finally this construction can be repeated for all $x \in Y$.  
For each $x\in Y$ there is a $x \in S_x \subseteq Y$ and there is a
map 
\[
\psi_x: S_x \to X \text{ such that }
\psi_x = (\varphi|_{\im(\psi_x)})^{-1}.
\]
Since $\varphi$ is a bijection, all local inverses must agree.
That is, $\{\psi_x: S_x \to X\}$ is a matching family on the open
cover $\{S_x\}$ of $Y$.
Since the topology is subcanonical, there is a unique amalgamation
$\psi: Y \to X$.
It remains to show that $\psi$ is an inverse for $\varphi$.
\[
\text{For all }S_x, \quad \varphi \circ \psi|_{S_x} = \varphi \circ
\psi_x = \Id_{S_x}.
\]
Therefore $\varphi$ is an isomorphism in \C.
\end{proof}

\section{Model categories for local po-spaces}

\subsection{A model category for local po-spaces}

Using our results on $\LPS$,
Theorem~\ref{thmMain} will now follow directly from Jardine's model
structure (Theorem~\ref{thmJardine}).

\begin{proof}[Proof of Theorem~\ref{thmMain}]
The open dicovers induce a Grothendieck topology on the small category
$\LPS$. 
Applying Theorem~\ref{thmEnoughPoints}, the Grothendieck topos
$\Sh(\LPS)$ has enough points.
So by Jardine's Theorem (Theorem~\ref{thmJardine}), $\sPre(\LPS)$ has a proper,
simplicial, cellular model structure in which
\begin{itemize}
\item the cofibrations are the monomorphisms, i.e. the levelwise
  monomorphisms of presheaves, 
\item the weak equivalences are the stalkwise equivalences,
  and 
\item the fibrations are the morphisms which have the right lifting property
with respect to all trivial cofibrations.
\end{itemize}
Finally by Theorem~\ref{thmStalkwiseTopLPS} the weak equivalences
coming from \LPS (via the Yoneda embedding) are precisely the isomorphisms.
\end{proof}

\subsection{Localization}

Our main motivation for constructing a model category for local
po-spaces was to model concurrent systems. In particular we would like
to be able to define and understand equivalences of concurrent systems
using such a model category. However our model structure on
$\sPre(\LPS)$ does not have any non-trivial equivalences among the
morphisms coming from $\LPS$. To obtain a model category more directly
useful for studying concurrency, we need to localize with respect to a
set of morphisms. 
In particular we want morphisms which preserve certain computer-scientific
information. 

How to best choose such morphisms is an important question and has
been studied in~\cite{bubenikContextBRICS}. 
For the sake of simplicity that paper studied only the category \Pospc of
po-spaces (a subcategory of $\LPS$).
There it was shown that the set of morphisms which should be
equivalences depends on the \emph{context}. 
That is, instead of choosing equivalences for \Pospc one should be
choosing equivalences for the coslice category or undercategory
$\APospc$ of po-spaces under a po-space $A$, where $A$ is called the
context. 

This result can be easily extended to our setting.
First we remark that if we choose a local po-space $A$ then
the undercategory $\ALPS$ is the category whose objects are dimaps $\iota_M: A \to
(M,\bar{U})$ and whose morphisms are dimaps $f: (M,\bar{U}) \to
(N,\bar{V})$ such that the following diagram commutes:
\[
\xymatrix{
& A \ar[dl]_{\iota_{M}} \ar[dr]^{\iota_{N}} & \\
(M,\bar{U}) \ar[rr]^f & & (N,\bar{V})
}
\]
Next, $\bar{y}(A) \in \sPre(\LPS)$ and the undercategory $\yAsPreLPS$
is the category whose objects are morphisms of simplicial presheaves
$\iota_{\alpha}: \bar{y}(A) \to \alpha$ and whose morphisms are
morphisms of simplicial presheaves $f: \alpha \to \beta$ such that the
following diagram commutes:
\[
\xymatrix{
& \bar{y}(A) \ar[dl]_{\iota_{\alpha}} \ar[dr]^{\iota_{\beta}} & \\
\alpha \ar[rr]^f & & \beta
}
\]
Since $\bar{y}:\LPS \to \sPre(\LPS)$ is a functor 
\[
\bar{y}(\iota_M): \bar{y}(A) \to \bar{y}(M,\bar{U}) \text{ and }
\bar{y}(\iota_N) = \bar{y}(f \circ \iota_M) = \bar{y}(f) \circ
\bar{y}(\iota_M).
\]
Hence $\ALPS$ embeds as a subcategory of $\yAsPreLPS$.

Define morphisms in $\yAsPreLPS$ to be weak equivalences,
cofibrations and fibrations if and only if they are weak equivalence,
cofibrations and fibrations in $\sPre(\LPS)$.
Then this makes $\yAsPreLPS$ into a model category
(see~\cite[Theorem 7.6.5]{hirschhornBook}).

We will show that this model category is again proper and cellular.
We will need the following definitions and a theorem of Kan.

\begin{definition}
\begin{itemize}
\item
Let $\C$ be a category and $I$ be a set of maps in $\C.$
A \emph{relative $I$-cell complex} is a map that can be constructed by a
transfinite composition of pushouts of elements of $I$.
\item
An object $A \in \C$ is \emph{small relative to a collection of morphisms
$\mathcal{D}$} in $\C$ if there exists a cardinal $\kappa$ such that for all
regular cardinals $\lambda \geq \kappa$ and for all
$\lambda$-sequences
\[
X_0 \to X_1 \to X_2 \to \ldots \to X_{\beta} \to \ldots
\]
with $X_{\beta} \to X_{\beta+1}$ in $\mathcal{D}$ for $\beta+1 <
\lambda$,
the set map
\[
\colim_{\beta<\lambda} \C(A,X_{\beta}) \to \C( A,
\colim_{\beta<\lambda} X_{\beta})
\]
is an isomorphism.
\end{itemize}
\end{definition}

\begin{definition} \label{defnCofGen}
A model category $\M$ is cofibrantly generated if there are sets $I$ and
$J$ such that 
\begin{itemize}
\item the domains of $I$ are small relative to the relative $I$-cell complexes,
\item the domains of $J$ are small relative to the relative $J$-cell complexes,
\item the fibrations have the right lifting property with respect to
  $J$, and
\item the trivial fibrations have the right lifting property with
  respect to $I$.
\end{itemize}
We say that $\M$ is cofibrantly generated by $I$ and $J$.
\end{definition}

\begin{definition}
\begin{itemize}
\item
Let $\M$ be a model category cofibrantly generated by $I$ and $J$.
An object $A \in \M$ is \emph{compact} if there is a cardinal $\gamma$ such that
for all relative $I$-cell complexes $f:X \to Y$ with a particular
presentation, every map $A \to Y$ factors through a subcomplex of size
at most $\gamma$.
\item
$f:A\to B$ is an \emph{effective monomorphism} if $f$ is the equalizer
of the inclusions $B \rightrightarrows B \amalg_A B$.
\end{itemize}
\end{definition}

\begin{definition} \label{defnCellular}
A \emph{cellular} model category is a model category cofibrantly
generated by $I$ and $J$ such that 
\begin{itemize}
\item the domains and codomains of elements of $I$ and $J$ are compact,
\item the domains of elements of $J$ are small relative to relative
  $I$-cell complexes, and
\item the cofibrations are effective monomorphisms.
\end{itemize}
\end{definition}

\begin{theorem}[{\cite[Theorem 11.3.2]{hirschhornBook}}] \label{thmKan}
Let $\M$ be a model category cofibrantly generated by the sets $I$
and $J$, and let $\N$ be a bicomplete category such that there exists
a pair of adjoint functors $F: \M \adjn \N: U$.
Define $FI = \{Fu \ | \ u \in I\}$ and $FJ = \{Fv \ | \ v \in J\}$.
If 
\begin{enumerate}
\item
the domains of $FI$ and $FJ$ are small relative to $FI$-cell and
$FJ$-cell, respectively, and
\item
$U$ maps relative $FJ$-cell complexes to weak equivalences,
\end{enumerate}
then $\N$ has a model category structure cofibrantly generated by $FI$
and $FJ$ such that $f$ is a weak equivalence in $\N$ if and only if
$Uf$ is a weak equivalence in $\M$, and $(F,U)$ is a Quillen pair.
\end{theorem}

\begin{theorem}
  Let $\M$ be a model category and let $A \in \M$.
Then $\AM$ has a model structure where a morphism
$
\xymatrix@=3mm{
& A \ar[dl] \ar[dr] \\
B \ar[rr]^f & & C
}
$
is a weak equivalence, cofibration or fibration in $\AM$ if and only if
$f$ is a weak equivalence, cofibration or fibration, respectively, in
$\M$.
If $\M$ is proper, cofibrantly generated or cellular, then so is $\AM$.
\end{theorem}

\begin{remark}
  For a more detailed proof we invite the reader to regard Hirschhorn's
  note~\cite{hirschhorn:overcategories}.
\end{remark}

\begin{proof}
That $\AM$ has the stated model structure follows from the definitions
(see \cite[Theorem 7.6.5]{hirschhornBook}).

Pushouts and pullbacks in $\AM$ can be formed by taking pushouts and
pullbacks of the underlying morphisms in $\M$, and then taking the
induced maps from $A$.
It thus follows that if $\M$ is proper so is $\AM$.

Assume $\M$ is cofibrantly generated by $I$ and $J$.
The method for showing that $\AM$ is cofibrantly generated will be to apply
Theorem~\ref{thmKan} to the following adjoint functors:
\[ F: \M \adjn (\AM) :U 
\]
where for $B \in \M$ and $f:B \to C \in \M$,
\[
F(B) = 
\begin{aligned}
\xymatrix@=5mm{A \ar[d]^{i_1} \\ A \amalg B}
\end{aligned}
, \quad
F(f) = 
\begin{aligned}
\xymatrix@=5mm{& A \ar[dl]_{i_1} \ar[dr]^{i_1} \\
 A \amalg B \ar[rr]^{\Id \amalg f} & & A \amalg C}
\end{aligned}
\]
and $U$ is the forgetful functor
\[
U 
\left( 
\begin{aligned}
\xymatrix@=5mm{A \ar[d]^{\iota_B} \\ B}
\end{aligned}
\right)
\ = \ B, \quad
U
\left(
\begin{aligned}
\xymatrix@=5mm{& A \ar[dl]_{\iota_B} \ar[dr]^{\iota_C} \\
   B \ar[rr]^{f} & & C}
\end{aligned}
\right)
\ = \ B \xto{f} C. 
\]
Define $FI = \{Fu \ | \ u \in I\}$ and $FJ = \{Fv \ | \ v \in J\}$.

The main observation for the proof is that for a morphism $u$ in $\M$,
the pushout of $Fu$ is obtained from the pushout of $u$ in $\M$.
That is,
\[
\begin{aligned}
\xy
\xymatrix{
A \ar[drr] \ar[dr] \ar[ddr] \\
& A \amalg B \ar[d]^{\iota_X \amalg f} \ar[r]_{\Id \amalg u} & A
\amalg C \ar[d] \\
& X \ar[r] & P
}
\POS( 32 , -20.4);\POS( 32,-22.8 ) \connect@{-}
\POS( 32 , -20.4);\POS( 34.4,-20.4 ) \connect@{-}
\endxy
\end{aligned}
\text{ where $P$ is defined by }
\begin{aligned}
\xy
\xymatrix{B \ar[r]^u \ar[d]_f & C \ar[d] \\
X \ar[r] & P
}
\POS( 8.4 , -8);\POS( 8.4,-10.4 ) \connect@{-}
\POS( 8.4 , -8);\POS( 10.8,-8 ) \connect@{-}
\endxy
\end{aligned}
\]
From this it follows that for a set of morphisms $S$ in $\M$, the
underlying morphisms of a relative $FS$-complex are a relative
$S$-complex.

Hence the conditions on $\AM$ in Theorem~\ref{thmKan} and the definition
of a cellular model category (Definition~\ref{defnCellular}) are
all inherited from the corresponding conditions in $\M$.

Finally one can check that the model category structure given by
Theorem~\ref{thmKan} coincides with the one in the statement of the
theorem. 
\end{proof}

Let $\mathcal{M}$ denote the model structure above on $\yAsPreLPS$.
Since $\M$ is cellular we can apply left Bousfield
localization~\cite{hirschhornBook} to this model structure
$\mathcal{M}$ with respect to a set of morphisms which will preserve
the computer-scientific properties we are interested in. 
In~\cite{bubenikContextBRICS}, one inverted the set of
\emph{dihomotopy equivalences} in $\APospc$.
So in our setting we will let $I$ be the set of \emph{dihomotopy
  equivalences} in $\ALPS$ defined below.
We will invert the set $\mathcal{I} = \{\bar{y}(f) \ | \ f \in I\}
\subset \yAsPreLPS$.
\begin{definition}
\begin{itemize}
\item Let $\dI$ be the po-space $([0,1],\leq)$ where $\leq$ is the
  usual total order on $[0,1]$. 
Given dimaps $f,g:(M,\bar{U}) \to (N,\bar{V}) \in \ALPS$,
$\phi$ is a \emph{dihomotopy} from $f$ to $g$ if $\phi:(M,\bar{U})
\times \dI \to (N,\bar{V})$, $\phi|_{(M,\bar{U})\times \{0\}} = f$,
$\phi|_{(M,\bar{U})\times \{1\}} = g$,  and for all $a \in A$,
$\phi(\iota_M(a),t) = \iota_N(a)$.
In this case write $\phi: f \to g$.
\item
The symmetric, transitive closure of dihomotopy is an equivalence relation. 
Write $f \simeq g$ if there is a chain of dihomotopies $f \to f_1
\from f_2 \to \ldots \from f_n \to g$.
\item
A dimap $f: (M,\bar{U}) \to (N,\bar{V})$ is a \emph{dihomotopy
  equivalence} if there is a dimap $g:(N,\bar{V}) \to (M,\bar{U})$
such that $g \circ f \simeq \Id_M$ and $f \circ g \simeq \Id_N$.
\end{itemize}
\end{definition}

The left Bousfield localization of $\mathcal{M}$ with respect to
$\mathcal{I}$ provides a model structure on \yAsPreLPS 
in which the weak equivalences are the $\mathcal{I}$-local
equivalences (see~\cite{hirschhornBook}), the cofibrations are the
cofibrations in $\mathcal{M}$ and the fibrations are morphisms which
have the right lifting property with respect to morphisms which are
both cofibrations and $\mathcal{I}$-local equivalences.

\begin{theorem}[Theorem~\ref{thmLocalized}]
Let $\mathcal{I} = \{ \bar{y}(f) \ | \ f $ is a  directed
    homotopy equivalence rel $A \}$.
The category \yAsPreLPS has a left proper, cellular
model structure in which  
\begin{itemize}
\item the cofibrations are the monomorphisms, 
\item the weak equivalences are the $\mathcal{I}$-local equivalences, and
\item the fibrations are those morphisms which have the right lifting
  property with respect to monomorphisms which are $\mathcal{I}$-local
  equivalences. 
\end{itemize}
\end{theorem}

We claim that this model category provides a good model for studying
concurrency. 
An analysis of this model category will be the subject of future
research.

\appendix

\section{Hypercovers}


Suppose now $\mathbb{C}$ is small and equipped with a Grothendieck topology,
i.e. we have a site $\left(\mathbb{C},\tau \right)$. The $\check{\textrm{C}}\textrm{ech}$
structure $\spre {\mathbb{C}}_{\check{c}\left(\tau \right)}$ is obtained
from the projective structure by homotopically localizing the comparison
morphisms given by the $\check{\textrm{C}}\textrm{ech}$ covers with
respect to $\tau $ or, up-to homotopy, from the injective structure
by localizing at the same set of morphisms. 

\begin{definition}
\label{def:cech-nerve}Let $U=\left\{ U_{i}\xrightarrow{{u_{i}}}X\right\} _{i\in I}\in J\left(X\right)$
be a cover. Let $i_{p}\in I$ for each $0\leq p\leq n$ and $U_{i_{0}\ldots i_{n}}$
be the \emph{wide pullback} of the $u_{i_{p}}$'s, i.e. the limiting
object of the diagram 
\begin{center} 
$
\xymatrix{
U_{i_0} \ar[drr]^{u_{i_0}} & 
\cdots & 
U_{i_p} \ar[d]^{u_{i_p}} & 
\cdots & 
U_{i_n} \ar[dll]_{u_{i_n}} 
\\
&& X &&
}
$
\end{center}
The $\check{C}ech$ nerve $\check{U}$ of $U$ is the simplicial
presheaf given by\[
\check{U}_{n}\, \deq \, \coprod _{i_{0},\ldots ,i_{n}\in I}y\left(U_{i_{0}\ldots i_{n}}\right)\]

\end{definition}
\begin{remark}
For any $n\in \mathbb{N}$, $X\in \mathbb{C}$ and $U\in J\left(X\right)$
there is a morphism\[
u_{i_{0}\cdots i_{n}}:\, U_{i_{0}\ldots i_{n}}\rightarrow X\]
and a diagram of presheaves
\begin{center} 
$
\xymatrix{
\check{U}_n \ar[rr]^{E_{U,X,n}} & &
y \left (X \right )
\\
y \left ( U_{i_0,\ldots,i_n} \right ) 
\ar[u]^{in_{y \left (U_{i_0,\ldots,i_n} \right )}}
\ar[urr]_{y \left ( u_{i_0,\ldots,i_n} \right )}
} 
$
\end{center}
where
$E_{U,X,n}$ is given by universal property. The $E_{U,X,n}$ assemble
to a morphism of simplicial presheaves\[
E_{U,X}:\, \check{U}\rightarrow \kappa _{y\left(X\right)}\]
\end{remark}

\begin{remark}
Given $U\in J\left(X\right)$ seen as a subcategory of the slice $\mathbb{C}/X$,
there is the evident functor \[
\begin{array}{rlcl}
 \delta _{U}: & U & \rightarrow  & \spre {\mathbb{C}}\\
  & u_{i} & \mapsto  & \kappa _{y\left(U_{i}\right)}\end{array}
\]

\end{remark}
\begin{proposition}
Localizing $\spre {\mathbb{C}}_{inj}$
at the sets 
\begin{enumerateroman}
\item $\left\{ E_{U,X}\, \mid \, X\in \mathbb{C},U\in J\left(X\right)\right\} $
;
\item $\left\{ hocolim\left(\delta _{U}\right)\rightarrow \kappa _{y\left(X\right)}\, \mid \, X\in \mathbb{C},U\in J\left(X\right)\right\} $;
\item $\left\{ \kappa \left(\iota _{U}\right)\, \mid \, X\in \mathbb{C},U\in J\left(X\right)\right\} $where,
given $X\in \mathbb{C}$ and $R$ a sieve on $X$, $\iota _{R}:\, R\hookrightarrow y\left(X\right)$
is the corresponding inclusion of presheaves;
\item $\left\{ \eta _{F}:\, F\rightarrow j\left(F\right)\, \mid \, F\in \spre {\mathbb{C}}\right\} $where
$j:\, \spre {\mathbb{C}}\rightarrow \spre {\mathbb{C}}$ is the objectwise
sheafification functor;
\end{enumerateroman}
yields the same model structure $\spre {\mathbb{C}}_{\check{c}\left(\tau \right)}$.
The same holds for the projective version.
\end{proposition}


Finally, there is a model structure $\spre {\mathbb{C}}_{hyp\left(\tau
  \right)}$ 
obtained from the projective structure by homotopically
localizing at the set of the comparison morphisms given by hypercovers
with respect to $\tau $. 
This model structure is Quillen equivalent to Jardine's model
  structure (Theorem~\ref{thmJardine}) on
  $\sSetCop$~\cite[Theorem 1.2]{duggerHollanderIsaksenHypercoversaSP}. 
As with
the $\check{\textrm{C}}\textrm{ech}$ structure, there is also an
injective version. Since $\check{\textrm{C}}\textrm{ech}$ covers
are particular hypercovers, there is the series of inclusions\[
\mathcal{W}_{prj}\subseteq \mathcal{W}_{\check{c}\left(\tau \right)}\subseteq \mathcal{W}_{hyp\left(\tau \right)}\]
and a similar series for the injective version. It is in general the
case that $\mathcal{W}_{\check{c}\left(\tau \right)}\subsetneqq \mathcal{W}_{hyp\left(\tau \right)}$,
yet equality holds in some important particular cases like the smooth
Nisnevitch site 
(c.f. \cite[Example A10]{duggerHollanderIsaksenHypercoversaSP}). 
It is an interesting question whether or not
  $\mathcal{W}_{\check{c}(\tau)} = \mathcal{W}_{hyp(\tau)}$ for local
  po-spaces.



\end{document}